# On the values of the functions

# $\frac{1}{\pi} Arg\, \zeta\left(\frac{1}{2}+in\right)$ and $\frac{1}{\pi} Arg\, \Gamma\left(\frac{1}{4}+\frac{in}{2}\right)$

by
Simon Plouffe
September 19, 2013

### Résumé

Une analyse de la fonction $\frac{1}{\pi} Arg\, \zeta(\frac{1}{2}+in)$ est présentée, cette analyse a permis de trouver une expression générale pour cette fonction qui n'utilise que des fonctions élémentaires comme la partie entière et fractionnaire. Cette formule générale permet de mettre en lumière une remarque générale qui avait été faite par Freeman Dyson à propos des valeurs de la fonction $\zeta$ reliées aux quasi-cristaux. On retrouve les mêmes valeurs pour $\frac{1}{\pi} Arg\, \Gamma(\frac{1}{4}+\frac{in}{2})$ qui sont très voisines de la fonction $\zeta$. De plus, les termes apparaissant dans ces formules générales obéissent à une règle simple avec n, à savoir, le nombre de 0 en binaire de n à partir de la droite. Ces approximations permettent également d'obtenir une formule plus exacte qui est asymptotique. Tous les résultats ont été obtenus de façon empirique.

### Abstract

An analysis of the function $\frac{1}{\pi} Arg\, \zeta(\frac{1}{2}+in)$ is presented. This analysis permits to find a general expression for that function using elementary functions of floor and fractional part. These formulas bring light to a remark from Freeman Dyson which relates the values of the $\zeta$ function to quasi-crystals. We find these same values for another function which is very similar, namely $\frac{1}{\pi} Arg\, \Gamma(\frac{1}{4}+\frac{in}{2})$. These 2 sets of formula have a definite pattern, the n'th term is related to values like $\pi, \ln(\pi), \ln(2), …, \log(p)$, where p is a prime number. The coefficients are closed related to a certain sequence of numbers which counts the number of 0's from the right in the binary representation of n. These approximations are regular enough to deduce an asymptotic and precise formula. All results presented here are empirical.

On the values of the functions $\frac{1}{\pi} Arg\, \zeta\left(\frac{1}{2}+in\right)$ and $\frac{1}{\pi} Arg\, \Gamma\left(\frac{1}{4}+\frac{in}{2}\right)$

## Introduction

Riemann établissait une formule qui compte le nombre de zéros de la fonction $\zeta$ qui lorsque T tend vers l'infini est

$$N(T) = \frac{T}{2\pi}\log\frac{T}{2\pi e} + \frac{7}{8} + O(\log(T)). \qquad (1)$$

L'argument étant que si T est sur la droite ½ ou non, a été repris par von Mangoldt, Edwards et Hardy-Littlewood. Nous n'essayerons pas du tout d'en établir la preuve mais simplement de savoir combien de zéros y a-t-il dans chaque intervalle unité. En faisant la différentielle, en principe on obtient le nombre de zéros, on pose donc $T \to T+1$ et on fait la différence pour s'apercevoir que le terme d'erreur est très imprécis. Tellement qu'en fait pour savoir combien y a-t-il de zéros dans chaque intervalle, la formule est inutile.

Un certain Backlund en 1912 eut l'idée d'établir plus précisément que

$$N(T) = \frac{1}{\pi}\vartheta(T) + 1 + S(T), \qquad (2)$$

où $\vartheta(T)$ est la fonction de Riemann-Siegel et $S(T)$ est $\frac{1}{\pi} Arg\, \zeta(\frac{1}{2}+iT)$. Restait à trouver comment varie cette fonction d'angle avec la fonction $\zeta$. Mais il se trouve que cette fonction $S(T)$ varie énormément et a donné lieu à plusieurs idées de la part de Freeman Dyson, Alain Connes, Hugh Montgomery et M.V. Berry. (voir la bibliographie). Et plus récemment à André LeClair et Guilherme França de Cornell.

Définissons F(n), comme étant le nombre de zéros de la fonction $\zeta$ dans chaque intervalle unité. Ce qui nous donne la suite pour n ≥ 1,

$$F(n) = 0, 0, 0, 0, 0, 0, 0, 0, 0, 0, 0, 0, 0, 0, 1, 0, 0, 0, 0, 0, 0, 1, 0, 0, 0, 1, 0, \dots, \qquad (3)$$

la valeur est 2 à n=111, 150, 169, 223, et est 3 à n=5826, 5978, 6494. Il y a donc 2 zéros entre 111 et 112. On compte ici le nombre de valeurs dont la partie entière est 111. Comme on le sait le nombre de zéros est environ $\frac{2\pi n}{\log(n)}$. Ce qui équivaut à dire que pour chaque entier la valeur de F(n) devrait être de $\frac{\log(n)}{2\pi}$, toujours lorsque n→ ∞. La valeur étant très lentement croissante, il est aisé d'utiliser un dessin représentant les valeurs en tons de gris, 0 = blanc, 1= gris clair, etc. Étant donné qu'avec 1 image on peut représenter des milliards de valeurs.

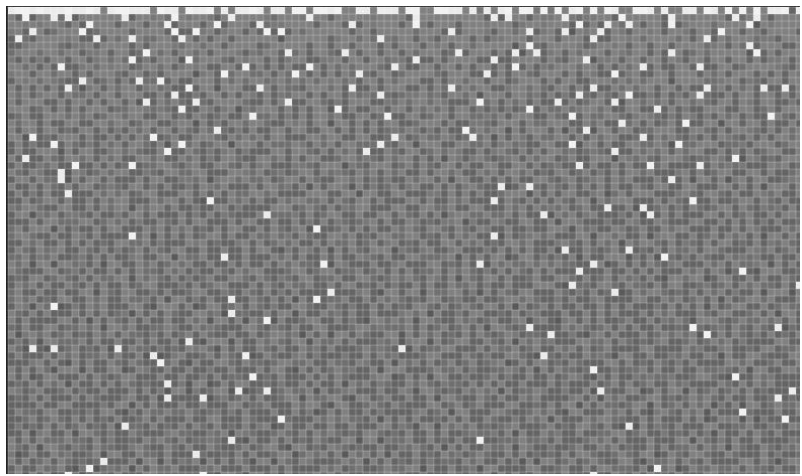

<p style="text-align:center">On the values of the functions $\frac{1}{\pi} Arg\, \zeta\left(\frac{1}{2}+in\right)$ and $\frac{1}{\pi} Arg\, \Gamma\left(\frac{1}{4}+\frac{in}{2}\right)$</p>

Ce dessin saute aux yeux, on y voit un motif immédiatement. Si on regarde attentivement le début sur la première ligne de pixels ont voit les zéros à 14, 21, 25, 30. La grille a été calibrée à 4000 pixels de largeur, le choix est arbitraire et le motif s'y trouve à toutes les valeurs. Si on fait le dessin aux valeurs de n>> 1 on trouve (colorisée).

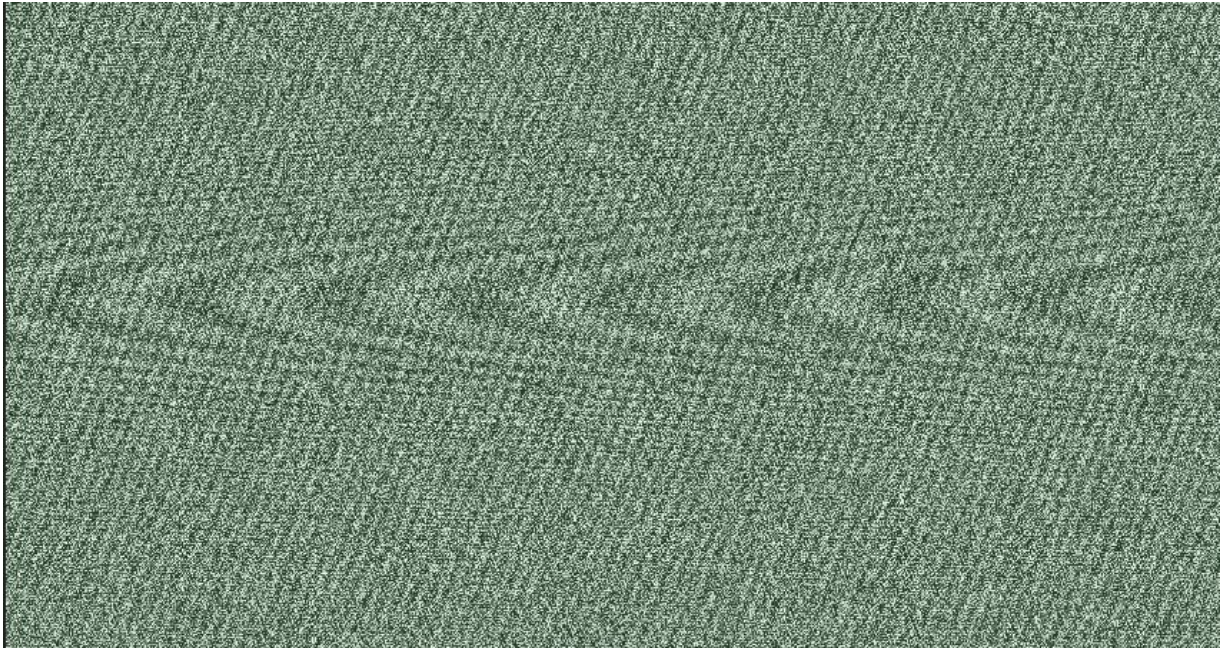

Le motif de vague est clair, c'est un Moiré. Les zéros de la fonction $\zeta$ sont des vagues. Ce n'est pas si surprenant, Granville l'a mentionné dans un article très intéressant et aussi Baillie. Les zéros peuvent être vus comme des ondes, en les sommant c'est une superposition d'ondes.

En variant le pas de l'image entre 4000 et 9064 ou 9065, on peut percevoir le premier battement, la valeur 9064 ou 9065 est en fait le premier terme d'une superposition de quelque chose dont on ignore la nature, on a $9064 < 1000 \frac{2\pi}{\log(2)} < 9065$. Mais quels sont les autres termes ? M. V. Berry a suggéré log(p) ou des puissances de log(p), p étant premier. Il suffit de rajouter ces termes mais comment ?

## Le nombre de zéros d'une fonction

Reprenons le calcul du nombre de zéros en commençant avec les fonctions les plus simples qui soit sont périodiques ou qui ont un grand nombre de zéros. La plus simple est sin(x), de façon évidente les zéros de sin(x) sont aux multiples de $\pi$, on prend ici x $\geq$ 0. Si on construit la fonction comme F(n) alors on aura : 0,0,1,0,0,1,0,0,1, …, en définissant

$$h(n,\alpha) = [(n+1)\alpha] - [n\alpha]. \qquad (4)$$

On a avec sin(x) que $\alpha = \frac{1}{\pi}$. Ce cas est trivial bien sûr mais nous servira plus loin. Prenons un autre exemple avec la fonction de Bessel $J_0(x)$. Les zéros se situent à intervalles d'à peu près $\pi$, le premier est 2.40483, 5.52008, 8.65373, 11.7915, 14.9309, … et pour ce qui nous intéresse, la suite 2,5,8,11,14, … Si on applique donc le même principe qu'avec sin(x) on trouve pas trop difficilement que le nombre de zéros B(n) de la fonction de Bessel $J_0(x)$ est donné par

On the values of the functions $\frac{1}{\pi} Arg\, \zeta\left(\frac{1}{2} + in\right)$ and $\frac{1}{\pi} Arg\, \Gamma\left(\frac{1}{4} + \frac{in}{2}\right)$

$$B(n) = h\left(n, \frac{4n}{\pi(4n-1)}\right) \tag{5}$$

Ce qui est conforme à ce qu'on connait de cette fonction dont les zéros sont espacés d'un multiple de $\pi$, petit à petit. Notons que la nature des zéros de cette fonction est inconnue. Ils ne semblent pas s'exprimer avec des constantes connues.

Prenons un exemple plus difficile comme la fonction de Airy, nous prendrons la valeur positive pour simplifier l'argument, appelons A(n) la fonction qui compte le nombre de zéros par intervalle unité. Alors,

$$A(n) = h\left(n, \frac{2n^{3/2}}{3\pi}\right) \tag{6}$$

Celle-ci a été obtenue en utilisant la formule asymptotique connue (voir DLMF) et en utilisant la méthode dite *bootstrap*. Même chose ici, on ne connait pas la nature des zéros de la fonction d'Airy mais on peut les compter par intervalle.

La question est maintenant, peut-on trouver une expression semblable avec $\zeta$ ?, si on utilise la formule asymptotique ça ne fonctionne pas, c'est trop imprécis. L'auteur a essayé une variété de méthodes *bootstrap* avec le log des premiers sans succès. Même si certaines pistes semblaient prometteuses. Si $y_n$ est la partie imaginaire du n'ième zéro non-trivial de la fonction $\zeta$ alors

$$y_1 \sim 2\frac{\pi}{(e^\pi - 1)},\ y_2 \sim 3\frac{\pi}{(e^\pi - 1)} \tag{7}$$

Mais là s'arrête le motif, ce n'est qu'une heureuse coïncidence.

## Densité ponctuelle

On appelle densité ponctuelle, le nombre de valeurs d'une fonction dans un intervalle donné, ici l'intervalle est l'unité. Pour obtenir la densité ponctuelle ou une fonction qui représente la densité ponctuelle d'une autre il est nécessaire d'inverser fonctionnellement.

Exemple, le nombre de nombres premiers est $\frac{n}{\log(n)}$ et inversement, la densité ponctuelle des nombres premiers est $\frac{1}{\log(n)}$. On peut bien sûr être plus précis comme le mentionne Dusart dans sa thèse. Pour la fonction $\zeta$, la densité est $\frac{2\pi n}{\log(n)}$, ce qui donne une densité ponctuelle de $\frac{\log(n)}{2\pi}$. Quand le terme principal est simple, la formule est directe mais si la formule est plus élaborée, c'est plus complexe. En reprenant la formule citée en [AL-GF] formule (13) on a le n'ième zéro y(n) de la fonction $\zeta$ est la solution de l'équation transcendante.

$$\frac{y(n)}{2\pi} \log \frac{y(n)}{2\pi e} + \frac{1}{\pi} Arg\, \zeta(\tfrac{1}{2} + iy(n)) = n - \frac{11}{8} \tag{8}$$

Si on résous, on obtient directement avec Maple ou Mathematica

On the values of the functions $\frac{1}{\pi} Arg\, \zeta\left(\frac{1}{2} + in\right)$ and $\frac{1}{\pi} Arg\, \Gamma\left(\frac{1}{4} + \frac{in}{2}\right)$

$$y(n) \approx 2\pi \frac{\left(n - \frac{11}{8}\right)}{W\left(\frac{\left(n - \frac{11}{8}\right)}{e}\right)}. \tag{9}$$

La fonction W est la dite 'de Lambert W', qui a des noms divers comme 'ProductLog'. La formule est étonnante et très précise quand n est grand, LeClair et França donnent des exemples ou n est de l'ordre de $10^{100}$ et $10^{200}$ avec une précision qui dépasse de loin tout ce qui a été fait avant. Ils notent également qu'ils ne tiennent pas compte des points de Gram ni non plus de la formule de Riemann-Siegel. Mais pour ce faire ils utilisent des méthodes numériques inversées, cette formule est le terme principal. La partie qui reste serait l'inverse fonctionnel du terme $\frac{1}{\pi} Arg\, \zeta(\frac{1}{2} + in)$. L'erreur avec la formule (9) seule est comprise entre -1 et 1. Les méthodes utilisées précédemment ne permettent pas de dégager une formule directe pour corriger le tir.

## Expressions pour $\frac{1}{\pi} Arg\, \zeta\left(\frac{1}{2} + in\right)$

Si on revient à la formule citée avec l'angle de la fonction $\frac{1}{\pi} Arg\, \zeta(\frac{1}{2} + in)$, en faisant le graphe on obtient, entre 1 et 115,

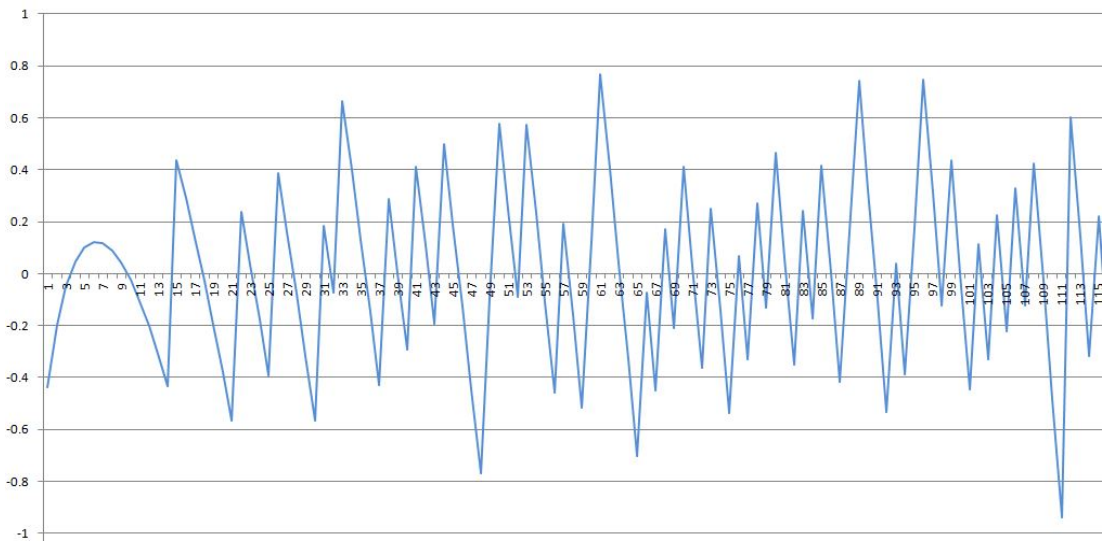

Les valeurs sont assez imprévisibles. De plus, l'argument de la somme de deux nombres complexes même simples se complique assez vite. Par exemple, en posant n=4000 et en faisant les sommes partielles de $\frac{1}{\pi} Arg\, \zeta(\frac{1}{2} + i4000)$, la somme à l'infini est -0.382343520341 mais le graphe des sommes partielles nous donne :

On the values of the functions $\frac{1}{\pi} Arg\, \zeta\left(\frac{1}{2}+in\right)$ and $\frac{1}{\pi} Arg\, \Gamma\left(\frac{1}{4}+\frac{in}{2}\right)$

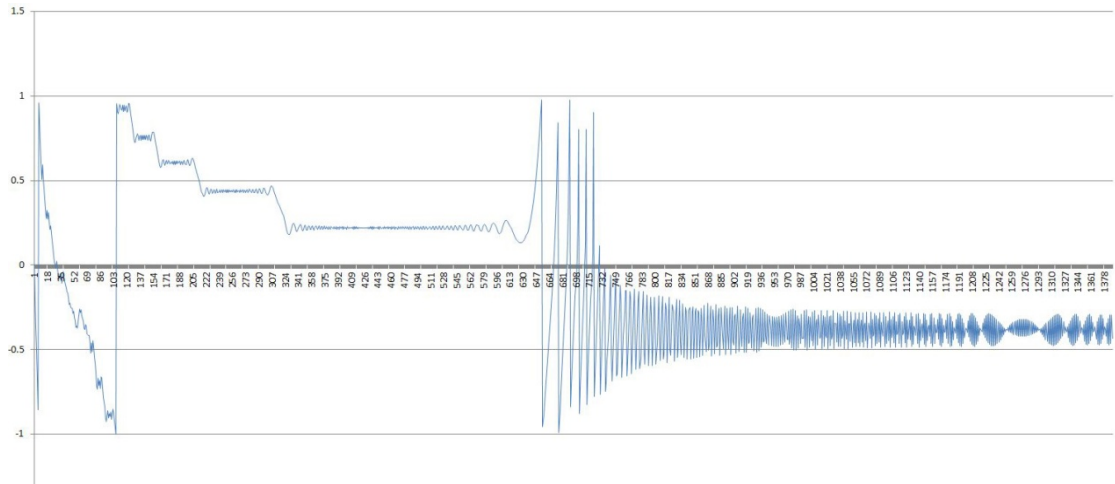

Donc, même si la somme ultimement est -0.382343520341, pour s'y rendre les valeurs bougent bien trop pour prédire quoi que ce soit. Cette courbe n'est pas caractéristique du comportement général, d'autres valeurs de n donnent des courbes complètement différentes.

On revient à la fonction qui donne approximativement le nombre de zéros,

$$N(T) \approx \frac{T}{2\pi} \log \frac{T}{2\pi e} + \frac{7}{8} \tag{10}$$

Si on rajoute le terme $\frac{1}{\pi} \text{Arg}\, \zeta(\frac{1}{2}+in)$, ou n remplace le T habituel, il apparaît évident que la somme des deux termes tend vers une valeur qui est toujours près de ½ en prenant la partie fractionnaire et ce de façon uniforme. Ce résultat est mentionné par LeClair et França. Mais, en regardant de près, la valeur ajoutée (l'angle avec ζ), si on pose que c'est exactement ½ alors ça implique que pour n=1 par exemple on a en (11) la valeur (qu'on appellera g(n)), une fois simplifiée nous donne,

$$g(1) = -\frac{1}{8\pi}(4\ln(2) + 4\ln(\pi) + 4 - 11\pi) = 0.92333784\ldots \tag{11}$$

La valeur de $\frac{1}{\pi} \text{Arg}\, \zeta(\frac{1}{2}+i)$ donne -0.43737200, la somme est 0.48596584 qui est proche de ½, donc la valeur approchée de $\frac{1}{\pi} \text{Arg}\, \zeta(\frac{1}{2}+i)$ serait $\frac{1}{8\pi}(4\ln(2) + 4\ln(\pi) + 4 - 7\pi)$.

$$g(2) = -\frac{1}{8\pi}(8\ln(\pi) + 8 - 11\pi) = 0.69231130\ldots \tag{12}$$

En appliquant le même principe on trouve donc la valeur approchée de $\frac{1}{\pi} \text{Arg}\, \zeta(\frac{1}{2}+2i) = -0.19231128$ qui comparée à la valeur exacte -0.19597758. Pour n = 10000 on trouve la valeur approchée de $\frac{1}{\pi} \text{Arg}\, \zeta(\frac{1}{2}+10000i)=$

$$-\frac{1}{8\pi}(120000\log(2) - 160000\log(5) + 40000\ln(\pi) + 40000 + 81129\,\pi) \tag{13}$$

On the values of the functions $\frac{1}{\pi} Arg\, \zeta\left(\frac{1}{2}+in\right)$ and $\frac{1}{\pi} Arg\, \Gamma\left(\frac{1}{4}+\frac{in}{2}\right)$

On remarque que les valeurs sont des combinaisons linéaires de π, log(π), log(2) et plus généralement log(p), avec p premier. On remarque également qu'à chaque valeur on se rapproche de la valeur supposée de ½. L'erreur à n=10000 est 0.0000066314559660306, cette valeur est $\frac{1}{2\pi}$ en fait, multiplié par un facteur de 24. Plus précisément, si e(10000) est l'erreur à n=10000, alors

$$e(10000) = 0.00000663145596030654980534783027\ldots \text{ et}$$
$$\frac{1}{2\pi}240000 = 0.00000663145596216230565703682347383\ldots \quad (14)$$

L'erreur est alors de l'ordre de $10^{-15}$ il reste suffisamment d'informations pour deviner le prochain terme qui est 7/120. Si on revient à la valeur trouvée en n=1, 2,3, on voit que cette correction est valide à partir de n=2, à n=10000000, l'erreur avec les 2 termes correctifs est de l'ordre de $10^{-29}$. Le 3ème terme du développement asymptotique est 93/7!, le 4ème est 127/896.

Finalement, en reprenant chaque coefficient trouvé on peut écrire,

$$\frac{1}{\pi}\text{Arg}\,\zeta\left(\frac{1}{2}+in\right) = \left\{\frac{n}{2\pi}\log\left(\frac{n}{2\pi e}\right) + \frac{7}{8}\right\} - \left(\frac{1}{\pi}\left(\frac{1}{48n} + \frac{7}{5760n^3} + \frac{31}{80640n^5} + \frac{127}{430080n^7} + \cdots\right)\right). \quad (15)$$

Et l'erreur est très petite par rapport à n, pour n=$10^4$ on a une erreur de l'ordre de $10^{-81}$.

On reconnaît ici le développement asymptotique de la formule de Riemann-Siegel, toujours avec la fonction { } qui représente la partie fractionnaire. On peut consulter les 2 suites A036282 et A114721 de la référence OEIS à ce sujet.

Voici un tableau des valeurs de $\frac{1}{\pi}\text{Arg}\,\zeta\left(\frac{1}{2}+in\right)$ pour n=1-19.

| $\frac{1}{\pi}Arg\,\zeta\left(\frac{1}{2}+in\right)$ | vraie valeur | approximation | valeur app. |
|---|---|---|---|
| 1 | -0.437372012317 | $\frac{1}{8\pi}(4\ln(2) + 4\ln(\pi) + 4 - 7\pi)$ | -0.423337836994 |
| 2 | -0.195977582921 | $\frac{1}{8\pi}(8\ln(\pi) + 8 - 7\pi)$ | -0.192311274139 |
| 3 | -0.046800452442 | $\frac{1}{8\pi}(-12\ln(3) + 12\ln(2) + 12\ln(\pi) + 12 - 7\pi)$ | -0.0445622398292 |
| 4 | 0.0474416622590 | $\frac{1}{8\pi}(-16\ln(2) + 16\ln(\pi) + 16 - 7\pi)$ | 0.0491062514140 |
| 5 | 0.101231367921 | $\frac{1}{8\pi}(-20\ln(5) + 20\ln(2) + 20\ln(\pi) + 20 - 7\pi)$ | 0.102560818219 |
| 6 | 0.122861669195 | $\frac{1}{8\pi}(-24\ln(3) + 24\ln(\pi) + 24 - 7\pi)$ | 0.123968719880 |
| 7 | 0.117778121706 | $\frac{7}{8\pi}(-4\ln(7) + 4\ln(2) + 4\ln(\pi) + 4 - \pi)$ | 0.118726607797 |
| 8 | 0.0898404109029 | $\frac{1}{8\pi}(-64\ln(2) + 32\ln(\pi) + 32 - 7\pi)$ | 0.090670102212 |
| 9 | 0.0419297327910 | $\frac{1}{8\pi}(-72\ln(3) + 36\ln(2) - 36\ln(\pi) + 36 - 7\pi)$ | 0.042667093960 |
| 10 | -0.0237198979997 | $\frac{1}{8\pi}(-40\ln(5) + 40\ln(\pi) + 40 - 7\pi)$ | -0.023056364340 |
| 11 | -.105325100472 | $\frac{1}{8\pi}(-44\ln(11) + 44\ln(2) - 44\ln(\pi) + 44 - 7\pi)$ | -0.104721949447 |

On the values of the functions $\frac{1}{\pi} Arg\, \zeta\left(\frac{1}{2} + in\right)$ and $\frac{1}{\pi} Arg\, \Gamma\left(\frac{1}{4} + \frac{in}{2}\right)$

| 12 | -.201429006842 | $\frac{1}{8\pi}(-48\ln(3) - 48\ln(2) + 48\ln(\pi) + 48 - 7\pi)$ | -0.200876161160 |
| --- | --- | --- | --- |
| 13 | -.310818966587 | $\frac{1}{8\pi}(-52\ln(13) + 52\ln(2) + 52\ln(\pi) + 42 - 7\pi)$ | -0.310308678190 |
| 14 | -.432469802098 | $\frac{7}{8\pi}(-8\ln(7) + 8\ln(\pi) + 8 - \pi)$ | -0.431995985476 |
| 15 | .434496598552 | $\frac{1}{8\pi}(-60\ln(3) - 60\ln(5) + 60\ln(2) + 60\ln(\pi) + 60 + \pi)$ | 0.434938810386 |
| 16 | 0.290840842657 | $\frac{1}{8\pi}(-192\ln(2) + 64\ln(\pi) + 64 + \pi)$ | 0.291255403206 |
| 17 | 0.137228161449 | $\frac{1}{8\pi}(-68\ln(17) + 68\ln(2) + 68\ln(\pi) + 68 + \pi)$ | 0.137618325910 |
| 18 | -0.0257546940666 | $\frac{1}{8\pi}(-144\ln(3) + 72\ln(\pi) + 72 + \pi)$ | -0.025386213458 |
| 19 | -0.197586328497 | $\frac{1}{8\pi}(-76\ln(19) + 76\ln(2) + 76\ln(\pi) + 76 + \pi)$ | -0.197237248073 |

# Expressions pour $\frac{1}{\pi} Arg\, \zeta\left(\frac{1}{2} + in\right)$, partie 2

Les coefficients apparaissant dans la 3ème colonne sont particuliers. On a toujours le terme en $8\pi$, vient ensuite le terme en $\pi$, $\log(\pi), \log(2)$ et $\log(p)$, les premiers p, sont ceux de la factorisation de n.

| coefficient | formule | remarque |
| --- | --- | --- |
| $\log(\pi)$ | 4n | |
| $\pi$ | $8\left[\frac{1}{8}\frac{4n\log(2) + 4n\log(\pi) - 4n\log(n) + 4n - 7\pi}{\pi}\right] - 7\pi$ | |
| $\log(2)$ | $a(n) = \begin{cases} 4n \text{ si } n \equiv 1 \bmod 4 \\ 0 \text{ si } n \equiv 2 \bmod 4 \\ 4n \text{ si } n \equiv 3 \bmod 4 \end{cases}$ | $\frac{-a(8n)}{32n} = A085068$ |
| $\log(3)$ | $a(n) = \begin{cases} 0 \text{ si } n \equiv 1 \text{ ou } 2 \bmod 3 \\ 0 \text{ si } n \equiv 2 \bmod 4 \end{cases}$ | $\frac{-a(3n)}{12n} = A051064$ |

Le coefficient de log(2) est surprenant, si on fait le graphe des valeurs on obtient,

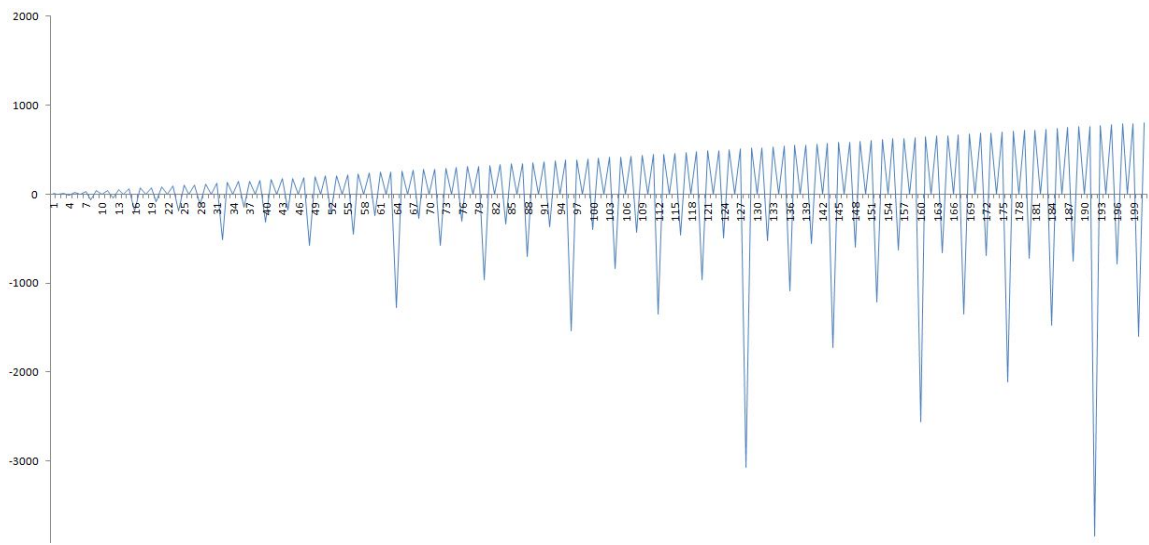

# On the values of the functions $\frac{1}{\pi} Arg\, \zeta\left(\frac{1}{2} + in\right)$ and $\frac{1}{\pi} Arg\, \Gamma\left(\frac{1}{4} + \frac{in}{2}\right)$

Les pics négatifs sont aux puissances de 2, les valeurs viennent en groupe de 3, suivis d'un zéro, la suite est 4, 0, 12, -16, 20, 0, 28, -64, 36, 0, 44, -48, 52, 0, 60, -192, 68, 0, 76, -80, 84, 0, 92, ... a(n) est parsemée de trous et a un maximum aux puissances de 2. En utilisant GFUN on trouve le lien assez facilement avec la suite A001511 ou A007814 du livre ou site de l'OEIS (voir bibliographie). La suite A001511 compte le nombre de zéros à partir de la droite en binaire de n. Si on prend $\frac{-a(8n)}{32n}$ on a la suite 2, 3, 2, 4, 2, 3, 2, 5, 2, 3, 2, 4, 2, 3, 2, 6, qui est A085058 du catalogue OEIS, c'est une variante de la suite générique A001511.

Le coefficient de log(3) suit la même règle en décalé, la suite donne : 0, 0, -12, 0, 0, -24, 0, 0, -72, 0, 0, -48, 0, 0, -60, 0, 0, -144, 0, 0, -84, 0, 0, -96, 0, 0, -324, 0, 0, -120, 0, 0, -132,... qui graphiquement donne,

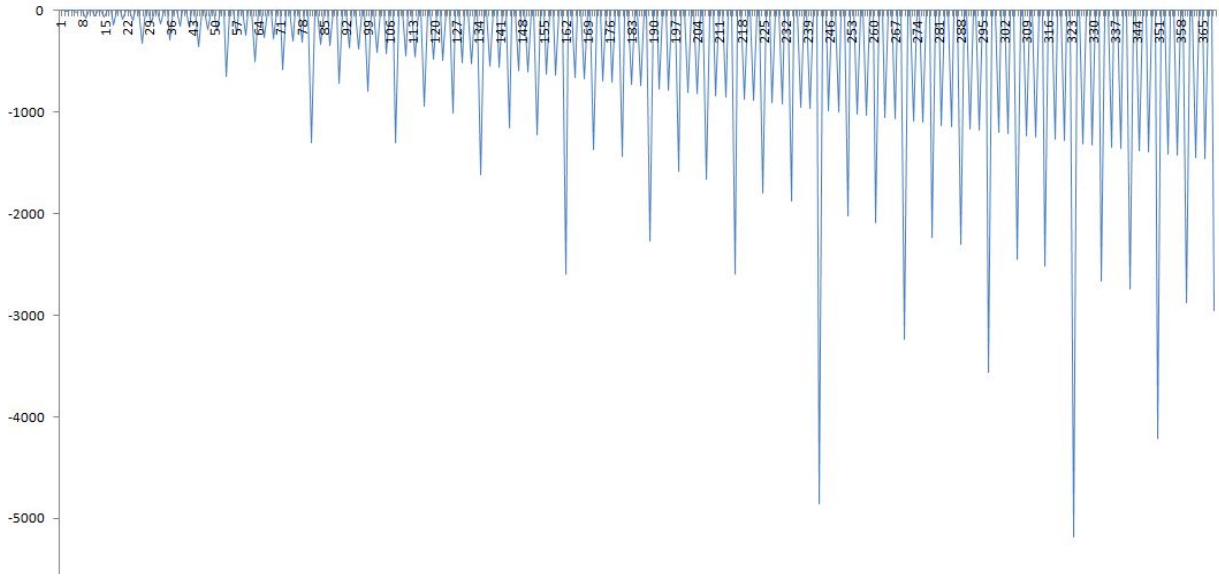

C'est le même motif fractal des suites A001511 ou A007814. Les autres termes en log(5), log(7) et plus généralement log(p) suivent le même motif exactement mais de plus en plus décalés. Si on prend un terme sur 3 et on pose $\frac{-a(3n)}{12n}$ inspiré par les résultats précédents on trouve : 1, 1, 2, 1, 1, 2, 1, 1, 3, 1, 1, 2, 1, 1, 2, 1, 1, 3, 1, 1 qui est la suite A051064 qui est de la même famille que la suite précédente. Les termes des autres rangs comme a(3n-1) et a(3n-2) sont triviaux ou 0. Le cas des coefficients de log(5), on retrouve la suite A091799 du catalogue OEIS.

On peut construire le même tableau pour la fonction $\frac{1}{\pi} Arg\, \Gamma\left(\frac{1}{4} + \frac{in}{2}\right)$ qui suit les mêmes règles. En fait, les 2 graphiques superposés de la fonction d'angle de $\zeta$ et de $\Gamma$ sont très proches, l'une complétant l'autre, dans [AL-GF], on donne la fonction équivalente à g(n) serait

$$g'(n) = \frac{n\log(\sqrt{\pi})}{\pi} \tag{16}$$

Ces 2 fonctions sont les *ondes porteuses* de ce phénomène + la partie avec l'argument de $\zeta$ et de $\Gamma$ qui suivent les règles énoncées plus haut.

On the values of the functions $\frac{1}{\pi} Arg\, \zeta\left(\frac{1}{2}+in\right)$ and $\frac{1}{\pi} Arg\, \Gamma\left(\frac{1}{4}+\frac{in}{2}\right)$

## Approximations de $\frac{1}{\pi} Arg\, \Gamma\left(\frac{1}{4}+\frac{in}{2}\right)$

| $\frac{1}{\pi} Arg\, \Gamma\left(\frac{1}{4}+\frac{in}{2}\right)$ | Valeur | Approximation | Valeur |
|---|---|---|---|
| 1 | -0.380438567847 | $\frac{-1}{8\pi}(4\ln(2)+\pi+4)$ | -0.394472743168 |

Étant donné que la fonction g(n) + $\frac{1}{\pi}\text{Arg}\,\zeta\left(\frac{1}{2}+in\right)$ compte le nombre de zéros, si on ne considère que la partie fractionnaire en $\frac{1}{\pi}\text{Arg}\,\zeta\left(\frac{1}{2}+in\right)$ c'est donc dire qu'à un multiple de π près ça compte aussi le nombre de zéros. Comme on le sait, $\frac{1}{\pi}\text{Arg}\,\zeta\left(\frac{1}{2}+in\right)$ varie entre -1 et 1. Donc, si effectue les différences, ∆ de cette fonction on obtient le décompte aussi, même si la valeur exacte comme tel est près d'un entier +1/2, il en est de même pour ∆., à partir de n=1, g(n) + $\frac{1}{\pi}\text{Arg}\,\zeta\left(\frac{1}{2}+in\right)$ =

0.485966, 0.496333, 0.497764, 0.498338, 0.498669, 0.498892, 0.499051, 0.499168, 0.499263, 0.499335, 0.499395, 0.499447, 0.499489, 0.499526, 1.49956, 1.49958, 1.49961, 1.49964, 1.49965, 1.49967, 1.49968, 2.49970, 2.49972, 2.49972, 2.49974, 3.49974, …

Qui arrondis et en faisant la différence donne bien le compte des zéros de la fonction F(n) originale 0, 0, 0, 0, 0, 0, 0, 0, 0, 0, 0, 0, 0, 1, 0, 0, 0, 0, 0, 0, 1, 0, 0, 0, 1, … qui comme le mentionne [AL-GF] manque une valeur entre 1007 et 1008 parce que + $\frac{1}{\pi}\text{Arg}\,\zeta\left(\frac{1}{2}+in\right)$ fait 1 tour de plus. La fonction donne bien l'angle exact mais au multiple de π près. Il en est de même avec la valeur approchée.

Dans les 2 cas, ça décompte bien le nombre de zéros, si la valeur diffère, il y a un tour ou plusieurs qui sont effectués par l'une ou l'autre expression mais le fond reste le même. Ce qui a été montré ici est que le *grain* ou la variation mystérieuse de la fonction $\frac{1}{\pi}\text{Arg}\,\zeta\left(\frac{1}{2}+in\right)$ et $\frac{1}{\pi}\text{Arg}\,\Gamma\left(\frac{1}{4}+\frac{in}{2}\right)$ sont en fait des approximations de constantes comme π, log(π), log(2) et log(p) avec l'aide de la fonction [ ] et { }. Donc, pour les premières valeurs de n et pour beaucoup d'autres, on a bien le décompte des zéros mais, à mesure que n croît, la fonction $\frac{1}{\pi}\text{Arg}\,\zeta\left(\frac{1}{2}+in\right)$ fait de plus en plus de tours à chaque entier.

Ce qui rejoint assez exactement les propos de Freeman Dyson à propos de quasi-cristaux. Il suggérait qu'à la base, les zéros ou les valeurs de la fonction ζ sont des variations générées à l'aide de constantes jumelées avec ces 2 fonctions [ ] et { }, que ce que l'on voit est un reflet.

On the values of the functions $\frac{1}{\pi} Arg\, \zeta\left(\frac{1}{2} + in\right)$ and $\frac{1}{\pi} Arg\, \Gamma\left(\frac{1}{4} + \frac{in}{2}\right)$

## Bibliographie et sur internet

On the values of the functions $\frac{1}{\pi} Arg\, \zeta\left(\frac{1}{2} + in\right)$ and $\frac{1}{\pi} Arg\, \Gamma\left(\frac{1}{4} + \frac{in}{2}\right)$